\documentclass[11pt]{article}

\usepackage{graphicx, color}
\usepackage{amsmath,amsfonts,amssymb}
\def\forall{\hbox{for all}~}
\def\L{\mathbf{L}}

\def\ve{\varepsilon}

\def\R{{\mathbb R}}

\def\vp{\varphi}

\def\v{\vskip 1em}

\def\C{{\cal C}}

\def\bega{\begin{array}}
\def\enda{\end{array}}
\def\begi{\begin{itemize}}
\def\endi{\end{itemize}}
\def\TV{\hbox{Tot.Var.}}

\def\ds{\displaystyle}
\def\bel{\begin{equation}\label}
\def\eeq{\end{equation}}
\def\sqr#1#2{\vbox{\hrule height .#2pt
\hbox{\vrule width .#2pt height #1pt \kern #1pt
\vrule width .#2pt}\hrule height .#2pt }}
\def\square{\sqr74}
\def\endproof{\hphantom{MM}\hfill\llap{$\square$}\goodbreak}

\begin{document}
\title{\bf Entropy Admissibility of the Limit Solution\\
for a Nonlocal Model of Traffic Flow}

\author{Alberto Bressan and Wen Shen\\ \, \\
Department of Mathematics, Penn State University.\\
University Park, PA~16802, USA.\\
\\
e-mails:~axb62@psu.edu,  
~wxs27@psu.edu}
\maketitle

\begin{abstract}
We consider a conservation law model of traffic flow, where the velocity of each car depends on a weighted average of the traffic density $\rho$ ahead. The averaging kernel is of exponential type: $w_\varepsilon(s)=\varepsilon^{-1} e^{-s/\varepsilon}$. 
For any  decreasing velocity function  $v$,  we prove that, as $\ve\to 0$, 
the limit of solutions to the nonlocal equation 
coincides with the unique entropy-admissible solution
to the scalar conservation law $\rho_t + (\rho v(\rho))_x=0$.
\end{abstract}

\section{Introduction}
\label{sec:1}
\setcounter{equation}{0}

We consider a nonlocal PDE model for traffic flow, where
the traffic density $\rho=\rho(t,x)$ 
satisfies a scalar conservation law with nonlocal flux
\bel{1} \rho_t+(\rho v(q))_x~=~0.\eeq
Here 
$\rho\mapsto v(\rho)$ is a decreasing function, modeling the velocity of cars
depending on the traffic density, while the integral
\bel{2}q(x)~=~\int_x^{+\infty} \ve^{-1} e^{(x-y)/\ve}\, \rho(y)\,ds\eeq
computes a weighted average of the density of cars ahead.
As in \cite{BS2},   we shall assume
\begi
\item[{\bf (A1)}] {\it The velocity function $v:[0,\rho_{jam}]\mapsto \R_+$ is $\C^2$,
and satisfies}
\bel{a1}
v(\rho_{jam}) \,=\, 0,\qquad\qquad v'(\rho)~\leq~-\delta_*~<~0, \qquad\forall \rho\in [0, \rho_{jam}].
\eeq\endi
One can think of  $\rho_{jam}$ as the maximum possible density of cars along the road, when all cars are packed bumper-to-bumper and nobody moves.  
The conservation equation (\ref{1}) will be solved with initial data
\bel{id}
\rho(0,x)~=~\bar \rho(x)~\in~[0,\rho_{jam}]\,.\eeq

As $\ve\to 0+$, the weight  function 
$w_\ve(s)=\ve^{-1} e^{-s/\ve}$ converges to a Dirac mass at the origin, 
and the nonlocal equation (\ref{1})-(\ref{2})
formally converges to the  scalar conservation law
\bel{Claw}
\rho_t + \bigl(\rho v(\rho)\bigr)_x~=~0. 
\eeq
Assuming that the initial datum $\bar \rho$ has bounded total variation and 
takes uniformly positive values, 
the recent analysis in \cite{BS2} has 
established:
\begi
\item[(i)] For every $\ve>0$, the Cauchy problem with non-local flux (\ref{1}), 
(\ref{2}), (\ref{id}),
has a unique solution $\rho=\rho_\ve(t,x)$. 
Its total variation satisfies a  uniform bound
\bel{TV}\TV\{ \rho_\ve(t,\cdot)\}~\leq~M\eeq
where the constant $M$ 
is independent of $t,\ve$.
\item[(ii)] As $\ve\to 0$, by possibly taking a subsequence,
one obtains  the convergence $\rho_\ve\to \rho$ in $\L^1_{loc}$.
The limit function $\rho=\rho(t,x)$ provides a weak solution to the 
Cauchy problem (\ref{id})-(\ref{Claw}).
\endi

A major issue, which was not fully resolved in \cite{BS2},
is the entropy admissibility of the 
limit solution $\rho$.   Aim of the present note is to resolve this question
in the affirmative.   Namely, we prove:
\v
{\bf Theorem.} {\it  Let $v$ satisfy the assumptions {\bf (A1)}, and 
let $\rho_\ve$ be a sequence of solutions to the nonlocal Cauchy problem
(\ref{1}), (\ref{2}) and (\ref{id}), satisfying the uniform BV bounds (\ref{TV}).
Assume that, as $\ve\to 0$, we have the convergence 
$\rho_\ve\to\rho$ in $\L^1_{loc}$.\\
Then $\rho$ is the unique entropy admissible solution to the Cauchy problem
(\ref{id})-(\ref{Claw}). }
\v
The above result was proved in \cite{BS2} in the special case
where the velocity is affine: $v(\rho) = a - b\rho$.   The earlier proof was based
on the Hardy-Littlewood inequality.   In the next section we give a simpler proof,
valid for a general class of velocity functions $v$.

For a more general class of averaging kernels, assuming that the initial datum
$\bar \rho$ satisfies a one-sided Lipschitz condition, the convergence to the unique
entropy admissible solution was recently proved in \cite{CCMS}.    Our result
requires an exponential kernel, but it applies to any BV initial data. In particular, 
$\bar \rho$ can be piecewise constant.

For the general theory of conservation laws we refer to 
\cite{Bbook, Daf, HR}.   A brief review of literature on hyperbolic
conservation laws with nonlocal flux can be found in \cite{BS2}.

\section{Proof of the theorem}
\label{s:2}
\setcounter{equation}{0}

{\bf 1.}
According to \cite{DLOW, Panov94}, to prove uniqueness
it suffices to prove that the limit solution dissipates one single strictly convex
entropy.
We thus consider the entropy and entropy flux pair
\bel{eep}\eta(\rho)\,=\,{\rho^2\over 2}\,,\qquad 
\qquad\psi(\rho)~=~\int_0^\rho \bigl[ s v(s) + s^2 v'(s)\bigr]\, ds.
\eeq
For future use, we observe that (\ref{2}) implies
\bel{dr}
\rho~=~q - \ve q_x.\eeq
Moreover, we introduce the function
\bel{VW}
W(\rho)~\doteq~\int_0^\rho s^2 v'(s)\, ds\,.\eeq
The equation  \eqref{1} can now be written as
\[ \rho_t + (\rho v(\rho) )_x ~=~ \Big(\rho (v(\rho) - v(q))\Big)_x\,.
\]
Multiplying both sides by $\eta'(\rho)=\rho$, we obtain
\bel{ee5}
 \eta(\rho)_t +\psi(\rho)_x ~=~ \rho\Big(\rho (v(\rho) - v(q))\Big)_x\,.
 \eeq
 \v
 {\bf 2.}
Given a test function $\vp\in \C^1_c(\R)$,  $\vp\geq 0$, using (\ref{dr})
we  estimate the quantity 
\bel{JJ}\bega{rl}
J&\ds \doteq~2\int  \rho\Big(\rho (v(\rho) - v(q))\Big)_x\vp\, dx\\[4mm]
&\ds=~\int  (\rho^2)_x \bigl(v(\rho) - v(q)\bigr)\,\vp\, dx
+\int  2\rho^2 \bigl(v(\rho) - v(q)\bigr)_x\vp\, dx \\[4mm]
&=\ds~- \int \rho^2\bigl(v(\rho) - v(q)\bigr)\vp_x\, ds
+ \int  \rho^2 \bigl(v(\rho) - v(q)\bigr)_x\,\vp\, dx\\[4mm]
&\doteq\ds~J_1+J_2\,.
\enda
\eeq
Concerning the second integral,
using (\ref{dr}) we obtain
\bel{J22} \bega{rl}J_2&\ds=~\int\rho^2 v'(\rho) \,\rho_x\, \vp\,dx - \int \rho q\, v'(q) q_x\, \vp\,dx
+\int \rho \ve\, (q_x)^2\, v'(q) \, \vp\,dx\\[4mm]
&\doteq~J_{21}+ J_{22} + J_{23}\,.\enda
\eeq
Using (\ref{dr}) once again,  we now compute
\bel{J222}\bega{rl} J_{21}+J_{22}&=~\ds
\int\rho^2 v'(\rho) \,\rho_x\, \vp\,dx - \int q^2\, v'(q) q_x\, \vp\,dx + \int q \ve\, (q_x)^2 
v'(q)\, \vp\, dx
\\[4mm]
&\doteq~J_3+J_4+J_5\,.
\enda
\eeq
Since $\rho,q,\vp\geq 0$ while $v'\leq 0$,
from (\ref{J22}) and (\ref{J222}) we immediately see that 
\bel{J23} J_{23}~\leq ~0,\qquad\qquad J_5~\leq~0.\eeq
On the other hand, integrating by parts and recalling (\ref{VW}), we obtain
\bel{J34}\bega{rl} J_3+ J_4\ds& =\ds~\int  \bigl[W(\rho)\bigr]_x \,\vp\, dx-\int \bigl[W(q)\bigr]_x \vp\, dx\\[4mm]
&=~\ds - \int \bigl[ W(\rho) - W(q)\bigr]\, \vp_x\, dx\,.\enda
\eeq
\v
{\bf 3.} To conclude, consider a sequence of solutions $\rho_\ve$ to
(\ref{1})-(\ref{2}), (\ref{id}). Assume that, as $\ve\to 0$, we have the convergence
$\rho_\ve\to \rho$ in $\L^1_{loc}$.   Notice that this implies
$q_\ve\to \rho$ in $\L^1_{loc}$ as well.  Hence, the integrals
$J_1$ and $J_3+J_4$  both approach zero.
By the previous analysis,
$$\bega{l}\ds 2
\int\!\! \int\bigl\{ \eta(\rho_\ve)\vp_t + \psi(\rho_\ve)\vp_x\bigr\}\, dxdt\\[4mm]
\qquad\ds 
\geq ~\int\!\! \int \rho_\ve^2\bigl(v(\rho_\ve)-v(q_\ve)\bigr)\vp_x\, dxdt +
\int\!\! \int  \bigl[ W(\rho_\ve) - W(q_\ve)\bigr]\, \vp_x\, dxdt\,.
\enda$$
Letting $\ve\to 0$, since the right hand side converges to zero, we obtain
$$\int\!\! \int\bigl\{ \eta(\rho)\vp_t + \psi(\rho)\vp_x\bigr\}\, dxdt~\geq~0.$$
This proves that the limit solution $\rho$ is entropy admissible.
In particular, by \cite{DLOW, Panov94}, $\rho$ is the unique entropy weak
solution to the Cauchy problem (\ref{id})-(\ref{Claw}).
\endproof

\v
{\bf Acknowledgment.}
This research  was partially supported by NSF with  
grant  DMS-2006884, ``Singularities and error bounds for hyperbolic equations".
\v

\end{document}